\documentclass[11pt,a4paper]{article}

\usepackage[T1]{fontenc}
\usepackage[utf8]{inputenc}
\usepackage{lmodern}
\usepackage[margin=1in]{geometry}
\usepackage{amsmath,amssymb,amsfonts,amsthm,mathtools}
\usepackage{microtype}
\usepackage{enumitem}
\usepackage{xcolor}
\usepackage{tikz}
\usetikzlibrary{arrows.meta,calc,fit}
\usepackage{hyperref}

\hypersetup{
  colorlinks=true,
  linkcolor=blue!70!black,
  citecolor=blue!70!black,
  urlcolor=blue!70!black,
  pdftitle={Counterexamples to the Laborde\textendash Payan\textendash Xuong conjecture and the Directed Path Partition Conjecture},
  pdfauthor={Jiangdong Ai, Xiangjun Wang, Ruibo Zhao},
  pdfsubject={Stable sets meeting all longest paths in oriented graphs},
  pdfkeywords={longest path, stable set, oriented graph, path partition conjecture, Laborde\textendash Payan\textendash Xuong conjecture}
}

\theoremstyle{plain}
\newtheorem{theorem}{Theorem}[section]
\newtheorem{proposition}[theorem]{Proposition}
\newtheorem{lemma}[theorem]{Lemma}
\newtheorem{corollary}[theorem]{Corollary}

\newtheorem{conjecture}[theorem]{Conjecture}

\theoremstyle{definition}

\newtheorem{question}[theorem]{Question}

\theoremstyle{remark}
\newtheorem{remark}[theorem]{Remark}

\newcommand{\ZZ}{\mathbb{Z}}

\tikzset{
  >={Stealth[length=2mm,width=1.5mm]},
  vtx/.style={circle,draw,inner sep=0pt,minimum size=6.4mm,font=\footnotesize,fill=white},
  tvtx/.style={vtx,fill=black!15},
  xvtx/.style={vtx,line width=1.1pt},
  wvtx/.style={draw,rounded corners=3.2mm,inner sep=3pt,minimum height=6.4mm,font=\footnotesize,fill=white},
  arc/.style={->,line width=0.5pt},
  qarc/.style={->,line width=1.5pt},
  all/.style={->,line width=1.6pt,black!60},
  grp/.style={draw,dashed,rounded corners=2mm,inner sep=5pt}
}

\title{Counterexamples to the Laborde--Payan--Xuong Conjecture and the Directed Path Partition Conjecture}

\author{Jiangdong Ai\thanks{School of Mathematical Sciences and LPMC, Nankai University, Tianjin 300071, China. \texttt{jd@nankai.edu.cn}. Supported by the National Natural Science Foundation of China (grant nos.~12522117 and 12401456) and the Fundamental and Interdisciplinary Disciplines Breakthrough Plan of the Ministry of Education of China (JYB2025XDXM207).}
\and Xiangjun Wang\thanks{School of Mathematics and Statistics, Huazhong University of Science and Technology, Wuhan 430074, China. \texttt{xjwang@hust.edu.cn}.}
\and Ruibo Zhao\thanks{School of Mathematics and Statistics, Huazhong University of Science and Technology, Wuhan 430074, China. \texttt{zorrorzhao@qq.com}.}}

\date{}

\begin{document}

\maketitle

\begin{abstract}
Let $D$ be a digraph and let $\lambda(D)$ denote the number of vertices in a longest path of $D$. In 1982, Laborde, Payan and Xuong conjectured that every digraph $D$ has a stable set which meets every longest path of $D$. For $\lambda(D)\ge2$, this is the case $q=1$ of the Directed Path Partition Conjecture, which states that for every digraph $D$ and every integer $q$ with $1\le q\le \lambda(D)-1$, there exists a partition $(A,B)$ of $V(D)$ such that $\lambda(D[A])\le q$ and $\lambda(D[B])\le \lambda(D)-q$. In this paper, we disprove both conjectures, even for oriented graphs. For every integer $\ell\ge 28$, we construct an oriented graph $D_\ell$ with $\ell(\ell+3)$ vertices and $\lambda(D_\ell)=\ell$ such that no stable set of $D_\ell$ meets every longest path of $D_\ell$. The smallest graph in this family is $D_{28}$, with $868$ vertices.


\end{abstract}

\section{Introduction}\label{sec:intro}

All digraphs in this paper are finite and have no loops or parallel arcs, but they may contain $2$-cycles, that is, pairs of opposite arcs. An \emph{oriented graph} is a digraph with no $2$-cycle. All paths and cycles are directed. A path is a non-empty sequence of distinct vertices with an arc from each vertex to the next.
The \emph{order} of a path $P$, denoted by $|P|$, is the number of its vertices. When a path is written as a concatenation of pieces, we allow some of the pieces to be empty; an empty piece contains no vertex and has order $0$.
For a non-empty digraph $D$, let $\lambda(D)$ denote the maximum order of
a path in $D$, and put $\lambda(D)=0$ if $D$ is empty. A path of order
$\lambda(D)$ is called a \emph{longest path} of $D$. Two vertices of $D$ are \emph{adjacent} if there is an arc between them, and a set of vertices of $D$ is \emph{stable} if no two of its vertices are adjacent. We say that a set $S\subseteq V(D)$ \emph{meets} a path $P$ if $S\cap V(P)\neq\emptyset$. For a set $Y\subseteq V(D)$, we write $D[Y]$ for the subdigraph of $D$ induced by $Y$ and $D-Y$ for $D[V(D)\setminus Y]$. For a vertex $v$ we write $D-v$ instead of $D-\{v\}$. The \emph{converse} of $D$ is the digraph obtained from $D$ by reversing all arcs. For a positive integer $n$, let $[n]=\{1,\dots,n\}$. Terminology and notation not introduced here follow~\cite{BJG}.

The Gallai--Hasse--Roy--Vitaver theorem~\cite{Gallai,Hasse,Roy,Vitaver} states that the chromatic number of the underlying graph of a digraph $D$ is at most $\lambda(D)$. Equivalently, $V(D)$ can be partitioned into at most $\lambda(D)$ stable sets. If $D$ is non-empty and $S$ is a stable set of $D$ which meets every longest path of $D$, then $\lambda(D-S)\le\lambda(D)-1$. Hence, if every digraph had such a stable set, then the theorem would follow at once by induction on $\lambda(D)$. Laborde, Payan and Xuong posed the following conjecture in 1982 and published it in 1983~\cite{LPX}.

\begin{conjecture}\label{conj:LPX}
Every digraph has a stable set which meets every longest path.
\end{conjecture}

Conjecture~\ref{conj:LPX} holds for acyclic digraphs, since the vertices of in-degree zero of an acyclic digraph form a stable set which contains the first vertex of every longest path. More generally, it holds for every digraph $D$ with a kernel, that is, a stable set $N$ such that every vertex outside $N$ has an out-neighbour in $N$. Indeed, the last vertex of a longest path either belongs to $N$ or has an out-neighbour in $N$. Any such out-neighbour must already lie on the path, since otherwise it would extend the path~\cite{GGM}. Laborde, Payan and Xuong~\cite{LPX} proved the conjecture for symmetric digraphs (those in which every arc has its reverse), and Havet~\cite{Havet} proved it for digraphs with stability number at most $2$. Further partial results, mostly for generalisations of tournaments, can be found in~\cite{GG,GGM,GR,vAFS,WW}. The conjecture is also listed in the Open Problem Garden~\cite{OPG}. It remained open for more than forty years.

Conjecture~\ref{conj:LPX} is closely related to the Path Partition Conjecture. For undirected graphs, this conjecture originated in a discussion between Lov\'asz and Mih\'ok in 1981 and first appeared in print in~\cite{LPX}. Frick~\cite{Frick} surveys its history. Bondy~\cite[Conjecture~4.45]{Bondy} stated a directed analogue with equality in both bounds. Frick, van Aardt, Dlamini, Dunbar and Oellermann~\cite{FADDO} studied the following formulation under the name Directed Path Partition Conjecture (DPPC).

\begin{conjecture}\label{conj:DPPC}
For every digraph $D$ and every integer $q$ with $1\le q\le\lambda(D)-1$, there exists a partition $(A,B)$ of $V(D)$ such that $\lambda(D[A])\le q$ and $\lambda(D[B])\le\lambda(D)-q$.
\end{conjecture}

The restriction of Conjecture~\ref{conj:DPPC} to oriented graphs is known as the Oriented Path Partition Conjecture (OPPC)~\cite{FADDO,Frick}.
Suppose that $\lambda(D)\ge2$ and let $A\subseteq V(D)$. Then $\lambda(D[A])\le1$ if and only
if $A$ is stable, and $\lambda(D-A)\le\lambda(D)-1$ if and only if
$A$ meets every longest path of $D$. Hence the case $q=1$ of
Conjecture~\ref{conj:DPPC} is precisely Conjecture~\ref{conj:LPX},
which holds trivially when $\lambda(D)\le1$. The DPPC has been verified for several classes of digraphs~\cite{AGGZ,AG,BJNY,BGLL,FADDO,GH}; see also the survey of Sun~\cite{Sun}. However, as recorded in~\cite{Frick}, it had not been settled even for $q=1$.

In this paper, we show that Conjecture~\ref{conj:LPX} is false, even for oriented graphs.

\begin{theorem}\label{thm:main}
For every integer $\ell\ge28$, there exists an oriented graph $D_\ell$ with $\ell(\ell+3)$ vertices and $\lambda(D_\ell)=\ell$ such that no stable set of $D_\ell$ meets every longest path of $D_\ell$.
\end{theorem}

The smallest of these oriented graphs, $D_{28}$, has $868$ vertices. By the observation above, Theorem~\ref{thm:main} has the following consequence.

\begin{corollary}\label{cor:main}
Conjecture~\ref{conj:LPX} is false, even for oriented graphs. Moreover, for every integer $\ell\ge28$, the oriented graph $D_\ell$ has no partition $(A,B)$ of its vertex set such that $\lambda(D_\ell[A])\le1$ and $\lambda(D_\ell[B])\le\ell-1$. In particular, the DPPC and the OPPC are false.
\end{corollary}

The same examples disprove the equality version stated by Bondy~\cite[Conjecture~4.45]{Bondy}. They also disprove the conjecture of Ai, Gerke, Gutin and Zhou~\cite{AGGZ} that every digraph has property BNY, a property first studied by Bang-Jensen, Nielsen and Yeo~\cite{BJNY}. Writing $\lambda_k(D)$ for the maximum number of vertices covered by $k$ pairwise vertex-disjoint paths of $D$, this property requires, for every $q$ as in Conjecture~\ref{conj:DPPC}, a partition $(A,B)$ of $V(D)$ such that $\lambda(D[A])\le q$ and $\lambda_k(D[B])\le\lambda_k(D)-q$ for every integer $k$ with $1\le k\le|B|$. For $k=1$ the latter is the inequality $\lambda(D[B])\le\lambda(D)-q$ of Conjecture~\ref{conj:DPPC}. On the other hand, the Path Partition Conjecture for undirected graphs is equivalent to Conjecture~\ref{conj:DPPC} restricted to symmetric digraphs, and it is not affected by our examples.

Let us describe the idea of the construction. Suppose that a digraph $D$ contains a subdigraph $H$ and a vertex
$z\notin V(H)$ such that $z$ is adjacent to every vertex of $H$ except
for the vertices of a directed triangle with vertex set $T$. Suppose
also that $H-t$ contains a longest path of $D$ for every $t\in T$.
Then no stable set $S$ which meets every longest path of $D$ contains
$z$. Indeed, if $z\in S$, then $S\cap V(H)\subseteq T$. Since the
vertices of $T$ are pairwise adjacent, $S\cap V(H)\subseteq\{t\}$
for some $t\in T$. Thus a longest path of $D$ in $H-t$ avoids $S$. If, in addition, $D$ has a longest path each of whose vertices plays the role of $z$ for some subdigraph of $D$, then no stable set meets every longest path of $D$. The difficulty is to attach these subdigraphs to the vertices of a path without creating paths which are too long. We overcome it by means of a gadget with $14$ vertices in which every path of maximum order among those starting at one of three prescribed vertices is forced to end on the triangle. Note that a directed cycle of odd length must occur in every counterexample to Conjecture~\ref{conj:LPX}, because a digraph without directed cycles of odd length has a kernel by a theorem of Richardson~\cite{Richardson}. In fact, by~\cite[Corollary~2.5]{FADDO}, the subdigraph induced by the first vertices of the longest paths of any counterexample contains a directed cycle of odd length.

The paper is organised as follows. In Section~\ref{sec:gadget} we define the gadget and prove its two key properties. In Section~\ref{sec:blocks} we extend the gadget to larger oriented graphs, which we call blocks. In Section~\ref{sec:proof} we define $D_\ell$ and prove Theorem~\ref{thm:main} and Corollary~\ref{cor:main}. Section~\ref{sec:remarks} contains further remarks and open questions.

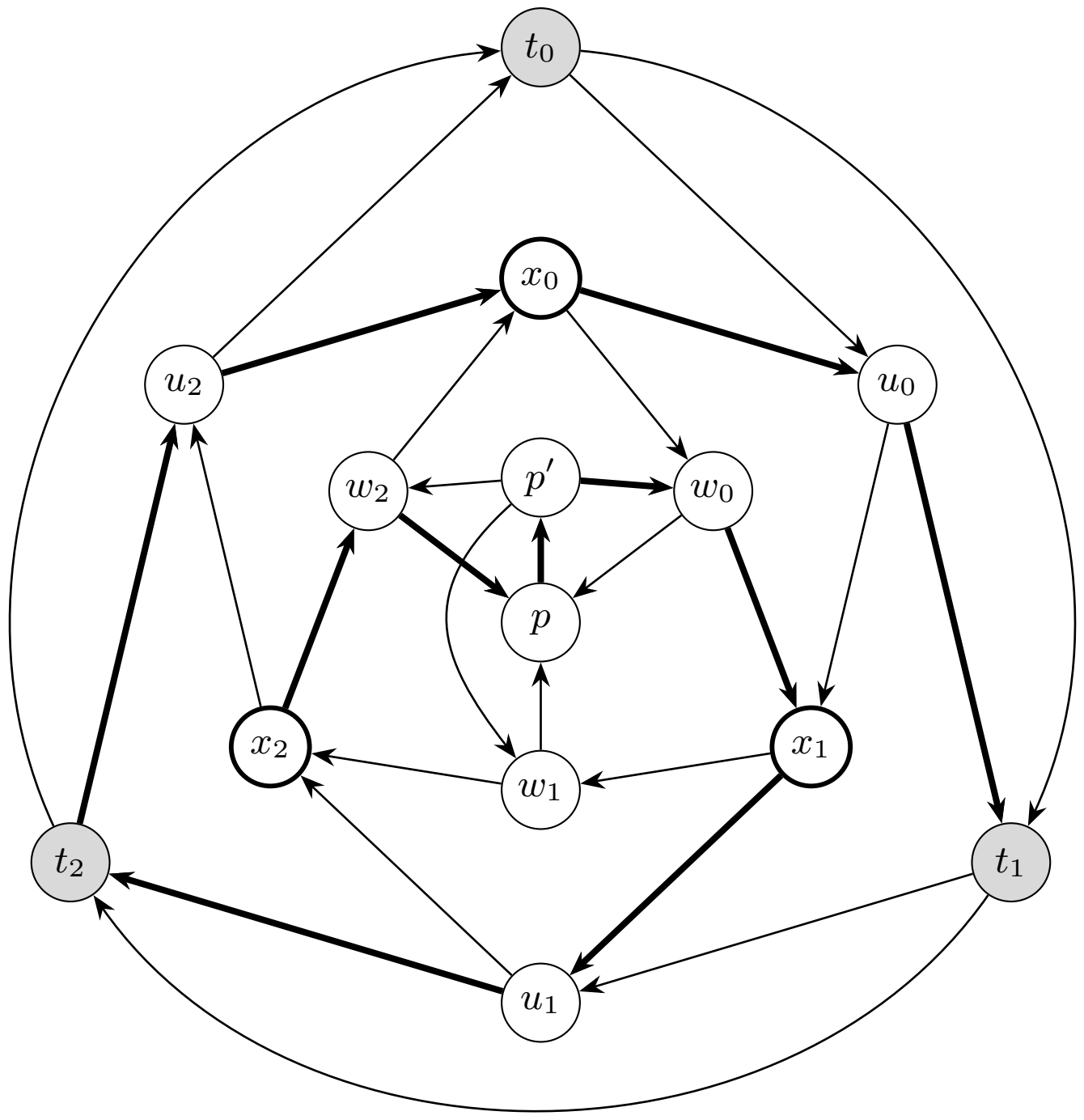
\begin{figure}[t]
\centering
\begin{tikzpicture}[scale=0.74]
\foreach \i/\a in {0/90,1/-30,2/-150}{
  \node[tvtx] (t\i) at (\a:6.0) {$t_{\i}$};
  \node[xvtx] (x\i) at (\a:3.45) {$x_{\i}$};
}
\foreach \i/\a in {0/30,1/-90,2/150}{
  \node[vtx] (u\i) at (\a:4.55) {$u_{\i}$};
  \node[vtx] (w\i) at (\a:2.2) {$w_{\i}$};
}
\node[vtx] (p) at (0,-0.35) {$p$};
\node[vtx] (q) at (0,1.25) {$p'$};
\foreach \i/\j in {0/1,1/2,2/0}{
  \draw[arc] (t\i) to[bend left=55] (t\j);
  \draw[arc] (t\i) -- (u\i);
  \draw[arc] (u\i) -- (t\j);
  \draw[arc] (u\i) -- (x\j);
  \draw[arc] (x\i) -- (u\i);
  \draw[arc] (x\i) -- (w\i);
  \draw[arc] (w\i) -- (x\j);
  \draw[arc] (w\i) -- (p);
}
\draw[arc] (p) -- (q);
\draw[arc] (q) -- (w0);
\draw[arc] (q) -- (w2);
\draw[arc] (q) to[out=222,in=128,looseness=1.3] (w1);
\draw[qarc] (x2) -- (w2);
\draw[qarc] (w2) -- (p);
\draw[qarc] (p) -- (q);
\draw[qarc] (q) -- (w0);
\draw[qarc] (w0) -- (x1);
\draw[qarc] (x1) -- (u1);
\draw[qarc] (u1) -- (t2);
\draw[qarc] (t2) -- (u2);
\draw[qarc] (u2) -- (x0);
\draw[qarc] (x0) -- (u0);
\draw[qarc] (u0) -- (t1);
\end{tikzpicture}
\caption{The gadget $G$. The vertices of $T$ are shaded and the vertices of $X$ are drawn with thick circles. The thick arcs form the path $Q_0=x_2w_2pp'w_0x_1u_1t_2u_2x_0u_0t_1$ of Lemma~\ref{lem:witness}, which has order $12$ and avoids $t_0$ and $w_1$. This diagram was drawn using OpenAI's GPT-6.2 Sol and checked by the authors.}
\label{fig:gadget}
\end{figure}

\section{The gadget}\label{sec:gadget}

Throughout the paper, the indices of the vertices $t_i$, $u_i$, $x_i$ and $w_i$ are taken modulo $3$. Let $G$ be the digraph with vertex set $T\cup U\cup X\cup W\cup\{p,p'\}$, where
\[
T=\{t_0,t_1,t_2\},\qquad U=\{u_0,u_1,u_2\},\qquad X=\{x_0,x_1,x_2\},\qquad W=\{w_0,w_1,w_2\},
\]
and with arc set
\[
A(G)=\{t_it_{i+1},\ t_iu_i,\ u_it_{i+1},\ u_ix_{i+1},\ x_iu_i,\ x_iw_i,\ w_ix_{i+1},\ w_ip,\ p'w_i : i\in\ZZ_3\}\cup\{pp'\}.
\]
Thus $G$ is an oriented graph with $14$ vertices and $28$ arcs, as shown in Figure~\ref{fig:gadget}. For convenience, we list the out-neighbourhood $N^+(v)$ and the in-neighbourhood $N^-(v)$ of every vertex $v$ of $G$:
\begin{center}
\begin{tabular}{c|cccccc}
$v$ & $t_i$ & $u_i$ & $x_i$ & $w_i$ & $p$ & $p'$\\ \hline
$N^+(v)$ & $t_{i+1},\,u_i$ & $t_{i+1},\,x_{i+1}$ & $u_i,\,w_i$ & $x_{i+1},\,p$ & $p'$ & $w_0,\,w_1,\,w_2$\\
$N^-(v)$ & $t_{i-1},\,u_{i-1}$ & $t_i,\,x_i$ & $u_{i-1},\,w_{i-1}$ & $x_i,\,p'$ & $w_0,\,w_1,\,w_2$ & $p$
\end{tabular}
\end{center}
In particular, $T$ induces a directed triangle and $X$ is a stable set. Moreover, if we put
\[
Y^+=T\cup U\qquad\text{and}\qquad Y^-=W\cup\{p,p'\},
\]
then there is no arc of $G$ between $Y^+$ and $Y^-$.

The first lemma is the core of the construction. It says that a path of $G$ which starts in $X$ has at most $12$ vertices, and that it can have $12$ vertices only if it ends in $T$.

\begin{lemma}\label{lem:gadget}
Let $P$ be a path of $G$ whose first vertex belongs to $X$. Then $|P|\le12$. Moreover, if $|P|=12$, then the last vertex of $P$ belongs to $T$.
\end{lemma}

\begin{proof}
Let $x^1,\dots,x^s$ be the vertices of $X$ on $P$, in the order in which they appear on $P$. Thus $1\le s\le3$ and $x^1$ is the first vertex of $P$. For $1\le r<s$, let $S_r$ be the subpath of $P$ consisting of the vertices strictly between $x^r$ and $x^{r+1}$, and let $S_s$ be the sequence of vertices after $x^s$ on $P$. The segment $S_s$ may be empty, while $S_1,\dots,S_{s-1}$ are non-empty because $X$ is stable. We call $S_1,\dots,S_s$ the \emph{segments} of $P$. Clearly,
\[
|P|=s+|S_1|+\dots+|S_s|.
\]
No segment contains a vertex of $X$, and there is no arc between $Y^+$ and $Y^-$. Hence each segment is contained either in $Y^+$ or in $Y^-$.

We first determine the segments $S_r$ with $r<s$. Let $r<s$ and let $x^r=x_i$. The first vertex of $S_r$ is an out-neighbour of $x_i$, so it is $u_i$ or $w_i$. The last vertex of $S_r$ is an in-neighbour of $x^{r+1}$. Since $x^{r+1}\ne x_i$ and $N^-(x_{j+1})=\{u_j,w_j\}$ for every $j\in\ZZ_3$, the last vertex of $S_r$ is $u_j$ or $w_j$ for some $j\in\{i,i+1\}$, and $x^{r+1}=x_{j+1}$. We call $S_r$ \emph{short} if $j=i$ and \emph{long} if $j=i+1$.

Suppose that $S_r\subseteq Y^+$. Then $S_r$ starts at $u_i$ and ends at $u_j$. If $j=i$, then $S_r=u_i$. If $j=i+1$, then the second vertex of $S_r$ is $t_{i+1}$, the only out-neighbour of $u_i$ in $Y^+$, and the penultimate vertex of $S_r$ is again $t_{i+1}$, the only in-neighbour of $u_{i+1}$ in $Y^+$. Hence $S_r=u_it_{i+1}u_{i+1}$.

Suppose that $S_r\subseteq Y^-$. Then $S_r$ starts at $w_i$ and ends at $w_j$. If $j=i$, then $S_r=w_i$. If $j=i+1$, then the second and third vertices of $S_r$ are $p$ and $p'$, because $p$ is the only out-neighbour of $w_i$ in $Y^-$ and $p'$ is the only out-neighbour of $p$, and the penultimate vertex of $S_r$ is again $p'$, the only in-neighbour of $w_{i+1}$ in $Y^-$. Hence $S_r=w_ipp'w_{i+1}$.

To summarise, the following holds for every $r<s$, where $x^r=x_i$.
\begin{enumerate}[label=\textup{(\roman*)}]
\item\label{it:short} If $S_r$ is short, then $|S_r|=1$ and $x^{r+1}=x_{i+1}$.
\item\label{it:long} If $S_r$ is long, then $S_r$ is one of the paths $u_it_{i+1}u_{i+1}$ and $w_ipp'w_{i+1}$, and $x^{r+1}=x_{i+2}$.
\end{enumerate}

Next we bound the order of the last segment $S_s$. If $S_s\subseteq Y^+$, then $|S_s|\le|Y^+|=6$. Suppose that $S_s\subseteq Y^-$ and that $S_s$ is non-empty, and let $x^s=x_i$. Then the first vertex of $S_s$ is $w_i$. Since $p$ is the only out-neighbour in $Y^-$ of each vertex of $W$, the only paths of $G[Y^-]$ starting at $w_i$ are $w_i$, $w_ip$, $w_ipp'$, $w_ipp'w_{i+1}$ and $w_ipp'w_{i+2}$. Hence $|S_s|\le4$.

We now distinguish four cases. Let $x^1=x_i$.

\emph{Case 1: $s=1$.} Then $|P|=1+|S_1|\le7$.

\emph{Case 2: $s\ge2$ and the segments $S_1,\dots,S_{s-1}$ are short.} Then $|P|\le s+(s-1)+6\le11$ by~\ref{it:short}.

\emph{Case 3: $s=2$ and $S_1$ is long.} Then $x^2=x_{i+2}$ by~\ref{it:long}. First let $S_1=u_it_{i+1}u_{i+1}$. Then $S_2$ is contained either in $Y^+\setminus V(S_1)$, a set of three vertices, or in $Y^-$. Hence $|S_2|\le4$ and $|P|\le2+3+4=9$. Now let $S_1=w_ipp'w_{i+1}$. If $S_2\subseteq Y^-$, then $V(S_2)\subseteq\{w_{i+2}\}$ and $|P|\le2+4+1=7$. If $S_2\subseteq Y^+$, then $|P|\le2+4+6=12$. Suppose that $|P|=12$. Then $S_2$ is a path of $G[Y^+]$ which contains all six vertices of $Y^+$, and its first vertex is $u_{i+2}$, the only out-neighbour of $x^2=x_{i+2}$ in $Y^+$. Recall that, for every $j\in\ZZ_3$, the only in-neighbour of $u_j$ in $Y^+$ is $t_j$ and the only out-neighbour of $u_j$ in $Y^+$ is $t_{j+1}$. Hence the second vertex of $S_2$ is $t_i$. The vertex $u_i$ lies on $S_2$ and is not its first vertex, so it is preceded on $S_2$ by $t_i$; thus $u_i$ is the third vertex of $S_2$, and the fourth vertex is $t_{i+1}$. In the same way, $u_{i+1}$ is preceded on $S_2$ by $t_{i+1}$, so $u_{i+1}$ is the fifth vertex of $S_2$. Therefore $S_2=u_{i+2}t_iu_it_{i+1}u_{i+1}t_{i+2}$, and the last vertex of $P$ is $t_{i+2}\in T$.

\emph{Case 4: $s=3$ and at least one of $S_1$, $S_2$ is long.} If exactly one of $S_1$ and $S_2$ were long, then~\ref{it:short} and~\ref{it:long} would give $x^3=x_{i+3}=x_i=x^1$, which is impossible. Hence both $S_1$ and $S_2$ are long, and so $x^2=x_{i+2}$ and $x^3=x_{i+4}=x_{i+1}$. By~\ref{it:long},
\[
S_1\in\{u_it_{i+1}u_{i+1},\ w_ipp'w_{i+1}\}\qquad\text{and}\qquad S_2\in\{u_{i+2}t_iu_i,\ w_{i+2}pp'w_i\}.
\]
The paths $S_1$ and $S_2$ are vertex-disjoint. Since $u_it_{i+1}u_{i+1}$ and $u_{i+2}t_iu_i$ share the vertex $u_i$, and $w_ipp'w_{i+1}$ and $w_{i+2}pp'w_i$ share the vertex $p$, only two possibilities remain.

First let $S_1=u_it_{i+1}u_{i+1}$ and $S_2=w_{i+2}pp'w_i$. The out-neighbours of $x^3=x_{i+1}$ are $u_{i+1}$, which lies on $S_1$, and $w_{i+1}$, and both out-neighbours $x_{i+2}$ and $p$ of $w_{i+1}$ occur before $x^3$ on $P$. Hence $|S_3|\le1$ and $|P|\le3+3+4+1=11$.

Now let $S_1=w_ipp'w_{i+1}$ and $S_2=u_{i+2}t_iu_i$. The out-neighbours of $x^3=x_{i+1}$ are $w_{i+1}$, which lies on $S_1$, and $u_{i+1}$. The out-neighbours of $u_{i+1}$ are $x_{i+2}$, which occurs before $x^3$ on $P$, and $t_{i+2}$. Finally, both out-neighbours $t_i$ and $u_{i+2}$ of $t_{i+2}$ lie on $S_2$. Hence $S_3$ is empty or one of the paths $u_{i+1}$ and $u_{i+1}t_{i+2}$. Therefore $|P|\le3+4+3+2=12$, and if $|P|=12$, then the last vertex of $P$ is $t_{i+2}\in T$.

In all cases $|P|\le12$, and if $|P|=12$, then the last vertex of $P$ belongs to $T$.
\end{proof}

The second lemma shows that the bound of Lemma~\ref{lem:gadget} is attained, and that it is attained by paths avoiding any prescribed vertex of $T$.

\begin{lemma}\label{lem:witness}
For every $j\in\ZZ_3$, the sequence
\[
Q_j=x_{j+2}\,w_{j+2}\,p\,p'\,w_j\,x_{j+1}\,u_{j+1}\,t_{j+2}\,u_{j+2}\,x_j\,u_j\,t_{j+1}
\]
is a path of $G$ of order $12$ whose first vertex belongs to $X$ and whose last vertex belongs to $T$, and $Q_j$ does not contain $t_j$.
\end{lemma}

\begin{proof}
The eleven pairs of consecutive vertices of $Q_j$ are
\[
x_{j+2}w_{j+2},\ \ w_{j+2}p,\ \ pp',\ \ p'w_j,\ \ w_jx_{j+1},\ \ x_{j+1}u_{j+1},\ \ u_{j+1}t_{j+2},\ \ t_{j+2}u_{j+2},\ \ u_{j+2}x_j,\ \ x_ju_j,\ \ u_jt_{j+1},
\]
and all of them are arcs of $G$. The vertices of $Q_j$ are pairwise distinct: they are the three vertices of $X$, the three vertices of $U$, the vertices $w_j$ and $w_{j+2}$, the vertices $t_{j+1}$ and $t_{j+2}$, and $p$ and $p'$. In particular, $Q_j$ contains neither $t_j$ nor $w_{j+1}$.
\end{proof}

\begin{figure}[t]
\centering
\begin{tikzpicture}[scale=0.92]
\node[vtx] (a1) at (0,0) {$a_1$};
\node[vtx] (a2) at (1.2,0) {$a_2$};
\node[font=\footnotesize] (ad) at (2.25,0) {$\cdots$};
\node[wvtx] (ak) at (3.55,0) {$a_{k-12}$};
\draw[rounded corners=3mm] (4.95,-1.8) rectangle (9.95,1.8);
\node[font=\small] at (7.45,-2.15) {$G$};
\foreach \i/\y in {0/1.05,1/0,2/-1.05}{
  \node[xvtx] (x\i) at (5.75,\y) {$x_{\i}$};
  \node[tvtx] (t\i) at (9.15,\y) {$t_{\i}$};
}
\node[font=\footnotesize,align=center] (mid) at (7.3,0) {$U\cup W$\\$\cup\,\{p,p'\}$};
\node[vtx] (b1) at (11.3,0) {$b_1$};
\node[font=\footnotesize] (bd) at (12.35,0) {$\cdots$};
\node[wvtx] (bk) at (13.6,0) {$b_{\ell-k}$};
\node[vtx] (z) at (7.45,3.9) {$z$};
\draw[arc] (a1) -- (a2);
\draw[arc] (a2) -- (ad);
\draw[arc] (ad) -- (ak);
\draw[arc] (b1) -- (bd);
\draw[arc] (bd) -- (bk);
\foreach \i in {0,1,2}{
  \draw[arc] (ak) -- (x\i);
  \draw[arc] (t\i) -- (b1);
}
\node[grp,fit=(a1)(ak)] (Lbox) {};
\node[grp,fit=(x0)(x2)(mid)] (Jbox) {};
\node[grp,fit=(b1)(bk)] (Rbox) {};
\node[font=\small] at (1.75,-0.85) {$L$};
\node[font=\small] at (12.45,-0.85) {$R$};
\draw[all] (Lbox.north) to[out=90,in=190] (z);
\draw[all] (Jbox.north) to[out=90,in=250] (z);
\draw[all] (z) to[out=-10,in=90] (Rbox.north);
\end{tikzpicture}
\caption{The oriented graph $F_k$. The arcs inside $G$ are not shown. A thick grey arrow from a dashed box to $z$ represents all arcs from the vertices in the box to $z$, and the thick grey arrow from $z$ to the dashed box around $R$ represents all arcs from $z$ to the vertices of $R$. Among the vertices of $F_k-z$, precisely those in $T$ are not adjacent to $z$. When $k=\ell$, the path $R$ and its incident arrows are absent. When $k=\ell-1$, $R$ consists of a single vertex. This diagram was drawn using OpenAI's GPT-6.2 Sol and checked by the authors.}
\label{fig:block}
\end{figure}
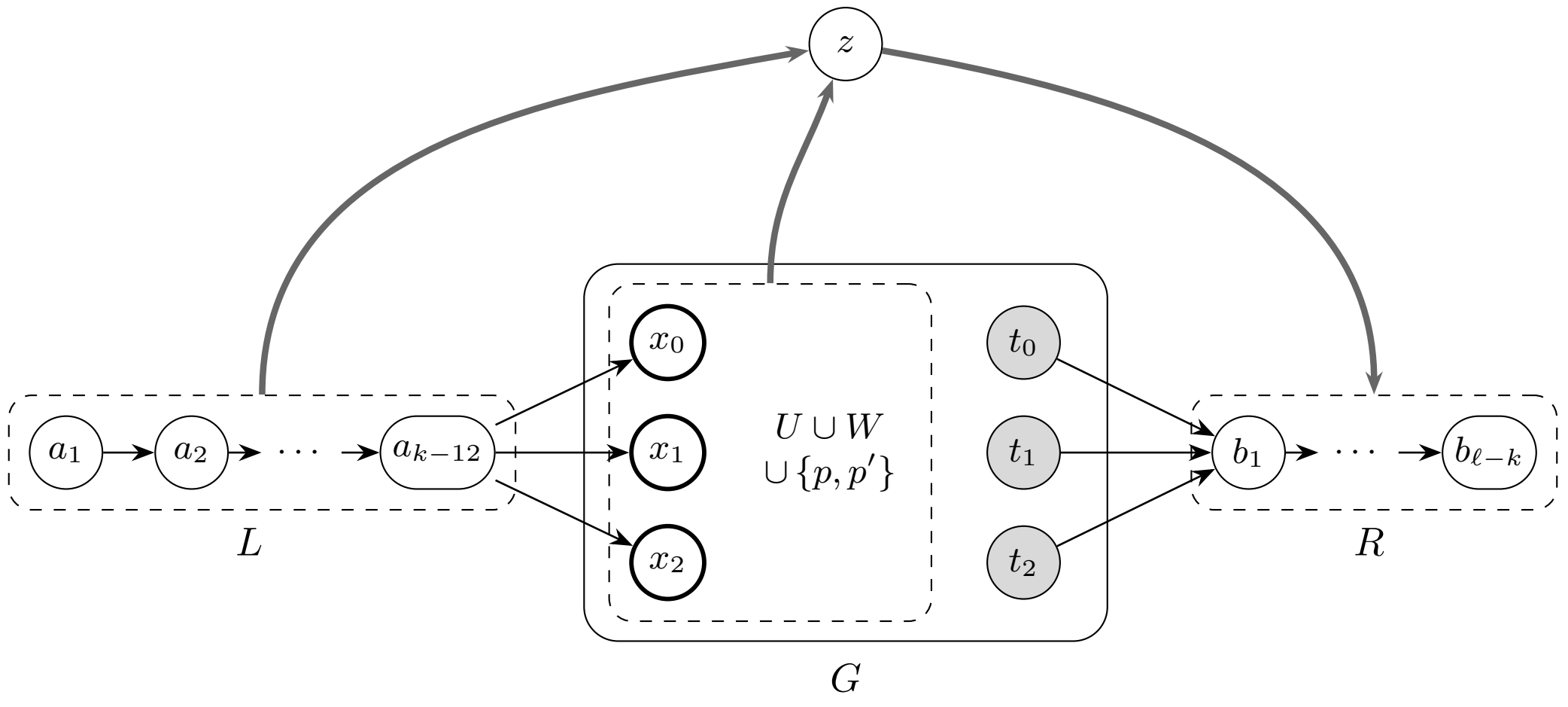

\section{The blocks}\label{sec:blocks}

Throughout Sections~\ref{sec:blocks} and~\ref{sec:proof}, fix an integer $\ell\ge28$. Let $k$ be an integer with $15\le k\le\ell$. We construct the oriented graph $F_k$ shown in Figure~\ref{fig:block}. Take a copy of $G$, a path $L=a_1a_2\dots a_{k-12}$ and a vertex $z$, all pairwise vertex-disjoint. If $k<\ell$, also take a path $R=b_1b_2\dots b_{\ell-k}$ disjoint from them. If $k=\ell$, let $R$ be the empty digraph. Then add the following arcs:
\begin{enumerate}[label=\textup{(F\arabic*)}]
\item\label{it:F1} $a_{k-12}x$ for every $x\in X$;
\item\label{it:F2} $tb_1$ for every $t\in T$, provided that $k<\ell$;
\item\label{it:F3} $vz$ for every $v\in V(L)\cup U\cup X\cup W\cup\{p,p'\}$;
\item\label{it:F4} $zb$ for every $b\in V(R)$.
\end{enumerate}
Note that $L$ has $k-12\ge3$ vertices and that $F_k$ has $14+(k-12)+(\ell-k)+1=\ell+3$ vertices. Since $G$ is an oriented graph, so is $F_k$. Let $H_k=F_k-z$.

The next lemma collects the properties of $F_k$ which we need. Parts~\ref{it:block-in} and~\ref{it:block-out} bound the order of a path of $H_k$ which can be extended through $z$, and parts~\ref{it:block-wit} and~\ref{it:block-adj} are the properties described in the introduction.

\begin{lemma}\label{lem:block}
Let $k$ be an integer with $15\le k\le\ell$. Then the following hold.
\begin{enumerate}[label=\textup{(\alph*)}]
\item\label{it:block-all} Every path of $H_k$ has order at most $\ell$.
\item\label{it:block-in} Every path of $H_k$ whose last vertex is an in-neighbour of $z$ in $F_k$ has order at most $k-1$.
\item\label{it:block-out} Every path of $H_k$ whose first vertex is an out-neighbour of $z$ in $F_k$ has order at most $\ell-k$.
\item\label{it:block-wit} For every $j\in\ZZ_3$, the oriented graph $H_k-t_j$ contains a path of order $\ell$.
\item\label{it:block-adj} The vertex $z$ is adjacent in $F_k$ to every vertex of $V(H_k)\setminus T$ and to no vertex of $T$.
\end{enumerate}
\end{lemma}

\begin{proof}
Part~\ref{it:block-adj} follows immediately from~\ref{it:F3} and~\ref{it:F4}.

In $H_k$, every arc with exactly one end in $V(L)$ is of the form $a_{k-12}x$ with $x\in X$, every arc with exactly one end in $V(R)$ is of the form $tb_1$ with $t\in T$, and there is no arc between $V(L)$ and $V(R)$. Thus no arc of $H_k$ enters $V(L)$ from outside and no arc of $H_k$ leaves $V(R)$. Consequently, on every path $P$ of $H_k$ the vertices of $L$ come first and the vertices of $R$ come last, so $P$ can be written as
\[
P=P_LP_GP_R,
\]
where $P_L$, $P_G$ and $P_R$ are the possibly empty portions of $P$ in $L$, $G$ and $R$, respectively. Each non-empty portion is a path in the corresponding digraph. If $P_L$ and $P_G$ are both non-empty, then the last vertex of $P_L$ is $a_{k-12}$ and the first vertex of $P_G$ belongs to $X$. If $P_G$ and $P_R$ are both non-empty, then the last vertex of $P_G$ belongs to $T$. If $P_G$ is empty, then $P_L$ or $P_R$ is empty.

\ref{it:block-all} Let $P=P_LP_GP_R$ be a path of $H_k$. If $P_L$ and $P_G$ are both non-empty, then the first vertex of $P_G$ belongs to $X$, so $|P_G|\le12$ by Lemma~\ref{lem:gadget}, and hence $|P|\le(k-12)+12+(\ell-k)=\ell$. If $P_L$ is empty, then $|P|\le14+(\ell-k)\le\ell-1$ because $k\ge15$. If $P_G$ is empty, then $|P|\le\max\{k-12,\ell-k\}<\ell$.

\ref{it:block-in} By~\ref{it:F3}, the in-neighbours of $z$ are the vertices of $V(L)\cup(V(G)\setminus T)$. Let $P=P_LP_GP_R$ be a path of $H_k$ whose last vertex is one of them. Then $P_R$ is empty. If $P_G$ is empty, then $|P|\le k-12$. If $P_L$ is empty, then $|P|\le14\le k-1$. If $P_L$ and $P_G$ are both non-empty, then the first vertex of $P_G$ belongs to $X$ and its last vertex does not belong to $T$. Hence $|P_G|\le11$ by Lemma~\ref{lem:gadget}, and $|P|\le(k-12)+11=k-1$.

\ref{it:block-out} By~\ref{it:F4}, the out-neighbours of $z$ are the vertices of $R$. A path of $H_k$ whose first vertex belongs to $R$ is a subpath of $R$, so its order is at most $\ell-k$.

\ref{it:block-wit} Let $j\in\ZZ_3$ and let $Q_j$ be the path of
Lemma~\ref{lem:witness}. By~\ref{it:F1}, the path $Q_j$ can be
preceded by $L$. If $k<\ell$, it can also be followed by $R$
by~\ref{it:F2}. The resulting path in $H_k$ avoids $t_j$ and has
order $(k-12)+12+(\ell-k)=\ell$.
\end{proof}

\section{Proof of Theorem~\ref{thm:main}}\label{sec:proof}

For $i\in[\ell]$, let
\[
k_i=\max\{i,\ell+1-i\}.
\]
Since $\ell\ge28$, we have $k_i\ge(\ell+1)/2>14$, and therefore $15\le k_i\le\ell$ for every $i\in[\ell]$. This is the only place where the assumption $\ell\ge28$ is used. Lemma~\ref{lem:block}\ref{it:block-in} would already fail at $k=14$
if we extended the definition of $F_k$ to that value. Indeed,
$t_0u_0t_1u_1t_2u_2x_0w_0x_1w_1x_2w_2pp'$ is a Hamiltonian path of $G$ ending at an in-neighbour of $z$. Note that $k_i=i$ if $2i\ge\ell+1$, and $k_i=\ell+1-i$ if $2i<\ell+1$.

\begin{figure}[t]
\centering
\begin{tikzpicture}[scale=0.95]
\node[vtx] (z1) at (0,0) {$z_1$};
\node[font=\footnotesize] (d1) at (1.5,0) {$\cdots$};
\node[vtx] (zi) at (3.3,0) {$z_i$};
\node[font=\footnotesize] (d2) at (5.5,0) {$\cdots$};
\node[vtx] (zj) at (8.4,0) {$z_j$};
\node[font=\footnotesize] (d3) at (10.9,0) {$\cdots$};
\node[vtx] (zl) at (12.4,0) {$z_\ell$};
\draw[arc] (z1) -- (d1);
\draw[arc] (d1) -- (zi);
\draw[arc] (zi) -- (d2);
\draw[arc] (d2) -- (zj);
\draw[arc] (zj) -- (d3);
\draw[arc] (d3) -- (zl);
\node[draw,rounded corners=2mm,font=\footnotesize,align=center,inner sep=4pt] (Ri) at (3.3,1.75) {copy of $R$\\$i-1$ vertices};
\node[draw,rounded corners=2mm,font=\footnotesize,align=center,inner sep=4pt] (Li) at (2.7,-1.8) {copies of $L$\\and of $G-T$};
\node[draw,rounded corners=2mm,font=\footnotesize,fill=black!15,inner sep=4pt,anchor=west] (Ti) at ($(Li.east)+(0.25,0)$) {$T^i$};
\node[draw,rounded corners=2mm,font=\footnotesize,align=center,inner sep=4pt] (Lj) at (7.8,1.8) {copies of $L$\\and of $G-T$};
\node[draw,rounded corners=2mm,font=\footnotesize,fill=black!15,inner sep=4pt,anchor=west] (Tj) at ($(Lj.east)+(0.25,0)$) {$T^j$};
\node[draw,rounded corners=2mm,font=\footnotesize,align=center,inner sep=4pt] (Rj) at (8.4,-1.75) {copy of $R$\\$\ell-j$ vertices};
\draw[all] (Ri) -- (zi);
\draw[all] (zi) -- (Li.north-|zi);
\draw[all] (Lj.south-|zj) -- (zj);
\draw[all] (zj) -- (Rj);
\node[font=\footnotesize] at (3.3,-2.95) {$F^i$ with $2i<\ell+1$};
\node[font=\footnotesize] at (8.4,-2.95) {$F^j$ with $2j\ge\ell+1$};
\end{tikzpicture}
\caption{Two blocks along the path $z_1z_2\dots z_\ell$. The left block $F^i$ is a copy of the converse of $F_{k_i}$ and the right block $F^j$ is a copy of $F_{k_j}$. Within $H^i$ and $H^j$, respectively, the unshaded box above each $z$-vertex contains its in-neighbours and the unshaded box below contains its out-neighbours. No vertex of either shaded triangle is adjacent to the $z$-vertex of its block. This diagram was drawn using OpenAI's GPT-6.2 Sol and checked by the authors.}
\label{fig:global}
\end{figure}
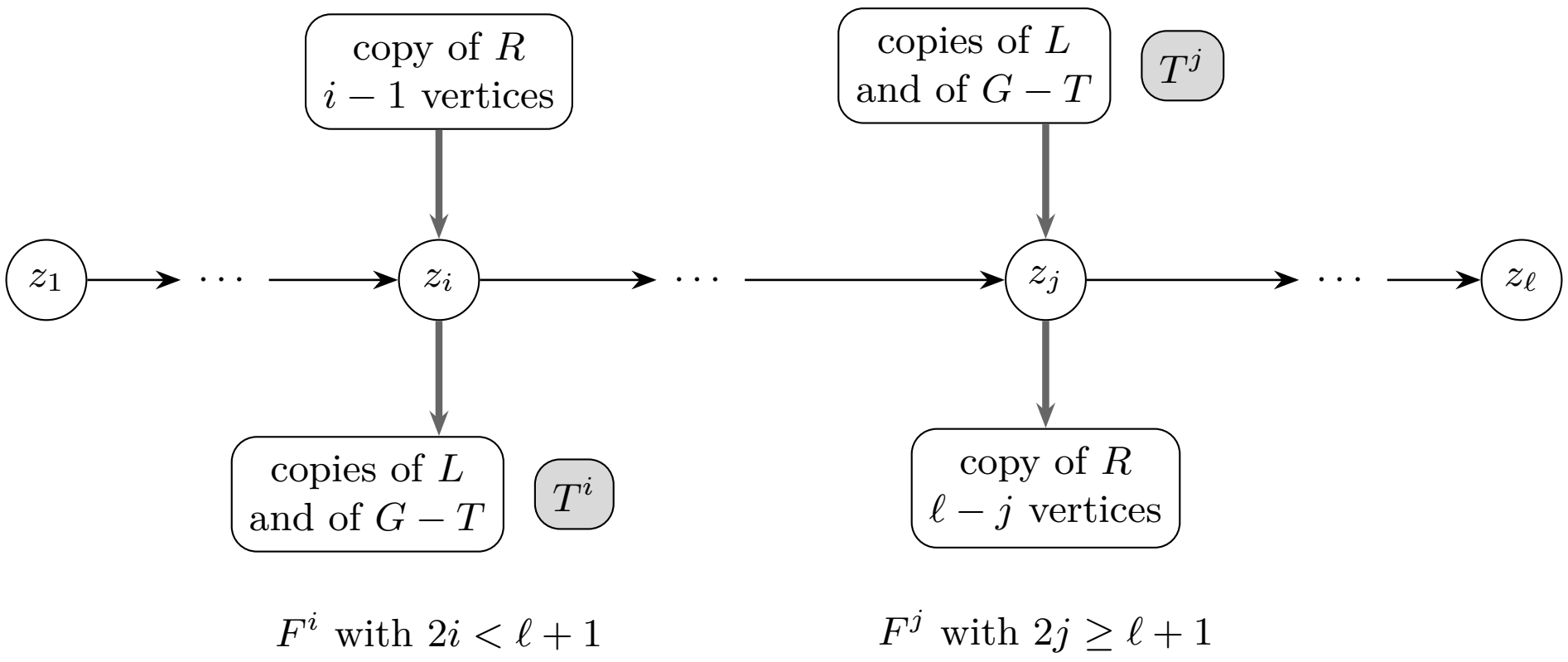

Let $F^1,\dots,F^\ell$ be pairwise vertex-disjoint oriented graphs such that, for every $i\in[\ell]$, the oriented graph $F^i$ is a copy of $F_{k_i}$ if $2i\ge\ell+1$, and a copy of the converse of $F_{k_i}$ if $2i<\ell+1$. Let $z_i$ and $T^i$ denote the copies in $F^i$ of the vertex $z$ and of the set $T$, respectively, and let $H^i=F^i-z_i$. We define $D_\ell$ to be the oriented graph obtained from the disjoint union of $F^1,\dots,F^\ell$ by adding the arcs $z_iz_{i+1}$ for all $i\in[\ell-1]$, as shown in Figure~\ref{fig:global}. Since each of $F^1,\dots,F^\ell$ has $\ell+3$ vertices, $D_\ell$ has $\ell(\ell+3)$ vertices. Let $Z=\{z_1,\dots,z_\ell\}$.

We first translate Lemma~\ref{lem:block} into properties of the oriented graphs $H^i$. The parameters $k_i$ and the orientations of the blocks are chosen precisely so that parts~\ref{it:Hi-in} and~\ref{it:Hi-out} of the following lemma hold.

\begin{lemma}\label{lem:Hi}
For every $i\in[\ell]$, the following hold.
\begin{enumerate}[label=\textup{(\alph*)}]
\item\label{it:Hi-all} Every path of $H^i$ has order at most $\ell$.
\item\label{it:Hi-in} Every path of $H^i$ whose last vertex is an in-neighbour of $z_i$ has order at most $i-1$.
\item\label{it:Hi-out} Every path of $H^i$ whose first vertex is an out-neighbour of $z_i$ has order at most $\ell-i$.
\item\label{it:Hi-wit} For every $t\in T^i$, the oriented graph $H^i-t$ contains a path of order $\ell$.
\item\label{it:Hi-adj} The vertex $z_i$ is adjacent to every vertex of $V(H^i)\setminus T^i$ and to no vertex of $T^i$.
\end{enumerate}
\end{lemma}

\begin{proof}
First let $2i\ge\ell+1$. Then $F^i$ is a copy of $F_{k_i}$ and $k_i=i$, so all five statements follow directly from Lemma~\ref{lem:block} with $k=i$.

Now let $2i<\ell+1$. Then $F^i$ is a copy of the converse of $F_{k_i}$ and $k_i=\ell+1-i$. The paths of $H^i$ are the reverses of the paths of a copy of $H_{k_i}$, the in-neighbours of $z_i$ in $F^i$ correspond to the out-neighbours of $z$ in $F_{k_i}$, and the out-neighbours of $z_i$ in $F^i$ correspond to the in-neighbours of $z$ in $F_{k_i}$. Hence~\ref{it:Hi-all}, \ref{it:Hi-wit} and~\ref{it:Hi-adj} follow from parts~\ref{it:block-all}, \ref{it:block-wit} and~\ref{it:block-adj} of Lemma~\ref{lem:block}. Moreover, a path of $H^i$ whose last vertex is an in-neighbour of $z_i$ is the reverse of a path of $H_{k_i}$ whose first vertex is an out-neighbour of $z$, so its order is at most $\ell-k_i=i-1$ by Lemma~\ref{lem:block}\ref{it:block-out}. Similarly, a path of $H^i$ whose first vertex is an out-neighbour of $z_i$ is the reverse of a path of $H_{k_i}$ whose last vertex is an in-neighbour of $z$, so its order is at most $k_i-1=\ell-i$ by Lemma~\ref{lem:block}\ref{it:block-in}.
\end{proof}

\begin{lemma}\label{lem:lambda}
We have $\lambda(D_\ell)=\ell$. Moreover, every path $P$ of $D_\ell$ which contains a vertex of $Z$ can be written as
\[
P=P_1\,z_iz_{i+1}\dots z_j\,P_2
\]
with $1\le i\le j\le\ell$, where $P_1$ is either empty or a path of $H^i$ whose last vertex is an in-neighbour of $z_i$, and $P_2$ is either empty or a path of $H^j$ whose first vertex is an out-neighbour of $z_j$. In particular, $|P_1|\le i-1$ and $|P_2|\le\ell-j$.
\end{lemma}

\begin{proof}
Since $z_1z_2\dots z_\ell$ is a path of $D_\ell$, we have $\lambda(D_\ell)\ge\ell$.

By the definition of $D_\ell$, every arc of $D_\ell$ whose ends lie in different oriented graphs among $F^1,\dots,F^\ell$ is of the form $z_iz_{i+1}$. Hence $D_\ell[Z]$ is the path $z_1z_2\dots z_\ell$, there is no arc between $V(H^i)$ and $V(H^r)$ for $i\ne r$, and $z_i$ is the only vertex outside $H^i$ which is adjacent to a vertex of $H^i$.

Let $P$ be a path of $D_\ell$. If $P$ contains no vertex of $Z$, then $P$ is a path of $H^i$ for some $i\in[\ell]$, and $|P|\le\ell$ by Lemma~\ref{lem:Hi}\ref{it:Hi-all}. Suppose that $P$ contains a vertex of $Z$. We claim that the vertices of $Z$ on $P$ are consecutive on $P$. Otherwise $P$ has a subpath $z'P'z''$, where $z'$ and $z''$ are distinct vertices of $Z$ and $P'$ is a non-empty path which contains no vertex of $Z$. Then $P'$ is a path of $H^r$ for some $r\in[\ell]$, and both $z'$ and $z''$ are adjacent to vertices of $H^r$, which is impossible. This proves the claim. Since $D_\ell[Z]$ is the path $z_1z_2\dots z_\ell$, these vertices form a path $z_iz_{i+1}\dots z_j$ with $i\le j$, and $P=P_1z_iz_{i+1}\dots z_jP_2$, where $P_1$ and $P_2$ contain no vertex of $Z$. If $P_1$ is non-empty, then its last vertex is an in-neighbour of $z_i$ outside $Z$, so $P_1$ is a path of $H^i$ and $|P_1|\le i-1$ by Lemma~\ref{lem:Hi}\ref{it:Hi-in}. If $P_2$ is non-empty, then its first vertex is an out-neighbour of $z_j$ outside $Z$, so $P_2$ is a path of $H^j$ and $|P_2|\le\ell-j$ by Lemma~\ref{lem:Hi}\ref{it:Hi-out}. Hence
\[
|P|=|P_1|+(j-i+1)+|P_2|\le(i-1)+(j-i+1)+(\ell-j)=\ell.
\]
Thus $\lambda(D_\ell)\le\ell$.
\end{proof}

\begin{proof}[Proof of Theorem~\ref{thm:main}]
Let $\ell\ge28$ and let $D_\ell$ be the oriented graph defined above. Then $D_\ell$ has $\ell(\ell+3)$ vertices, and $\lambda(D_\ell)=\ell$ by Lemma~\ref{lem:lambda}.

Suppose, for a contradiction, that $S$ is a stable set of $D_\ell$ which meets every longest path of $D_\ell$. Since $z_1z_2\dots z_\ell$ is a longest path of $D_\ell$, there is an index $i\in[\ell]$ such that $z_i\in S$. By Lemma~\ref{lem:Hi}\ref{it:Hi-adj}, the vertex $z_i$ is adjacent to every vertex of $V(H^i)\setminus T^i$, so $S\cap V(H^i)\subseteq T^i$. Since $T^i$ induces a directed triangle, $S$ contains at most one vertex of $T^i$. Hence there is a vertex $t\in T^i$ such that $S\cap V(H^i)\subseteq\{t\}$. By Lemma~\ref{lem:Hi}\ref{it:Hi-wit}, the oriented graph $H^i-t$ contains a path $Q$ of order $\ell$. Then $Q$ is a longest path of $D_\ell$ and $S\cap V(Q)=\emptyset$, a contradiction.
\end{proof}

\begin{proof}[Proof of Corollary~\ref{cor:main}]
Let $\ell\ge28$, and suppose that $(A,B)$ is a partition of $V(D_\ell)$ such that $\lambda(D_\ell[A])\le1$ and $\lambda(D_\ell[B])\le\ell-1$. Then $A$ is a stable set, and no longest path of $D_\ell$ is contained in $D_\ell[B]$, that is, $A$ meets every longest path of $D_\ell$. This contradicts Theorem~\ref{thm:main}. Hence the oriented graph $D_\ell$ is a counterexample to Conjecture~\ref{conj:LPX}, and to Conjecture~\ref{conj:DPPC} with $q=1$.
\end{proof}

\section{Further remarks}\label{sec:remarks}

\begin{remark}\label{rem:two-cycles}
If $2$-cycles are allowed, then slightly smaller examples can be
obtained. Let $G'$ be obtained from $G$ by identifying $p$ and $p'$
to form a single vertex $p^*$ and deleting the resulting loop. It has
$13$ vertices, and $p^*$ forms a $2$-cycle with each of $w_0$, $w_1$
and $w_2$.

Let $P$ be a path of $G'$ starting in $X$. Suppose first that $P$
contains $p^*$. The predecessor of $p^*$ on $P$, and its successor if
there is one, belong to $W$. Replacing $p^*$ by $pp'$ therefore gives
a path $\widehat{P}$ of $G$ of order $|P|+1$ starting in $X$. By
Lemma~\ref{lem:gadget}, $|P|+1\le12$, and if $|P|+1=12$, then the last
vertex of $\widehat{P}$ belongs to $T$. If $p^*$ is the last vertex of
$P$, then $\widehat{P}$ ends at $p'\notin T$, and so $|P|\le10$.
Otherwise $P$ and $\widehat{P}$ have the same last vertex. Hence
$|P|\le11$, and if $|P|=11$, then the last vertex of $P$ belongs to $T$.

Now suppose that $P$ avoids $p^*$. Then $P$ is a path of $G-\{p,p'\}$.
In the proof of Lemma~\ref{lem:gadget}, each non-empty segment in $Y^-$ now
consists of a single vertex of $W$. Case~1 still gives $|P|\le7$,
and in Case~3 the long segment lies in $Y^+$, so $|P|\le2+3+3=8$.
Case~4 is impossible because its two long segments would both lie
in $Y^+$ and share a vertex. In Case~2 we have $|P|\le11$, with
equality only if $s=3$ and the last segment uses all of $Y^+$.
Write $x^s=x_h$. The segment starts at $u_h$ and must be
$u_ht_{h+1}u_{h+1}t_{h+2}u_{h+2}t_h$, by the same in- and
out-neighbour restrictions as in Case~3. Thus equality again
forces the last vertex into $T$. Replacing $pp'$ by $p^*$ in the
paths $Q_j$ of Lemma~\ref{lem:witness} gives paths of $G'$ of order
$11$ which start in $X$, end in $T$ and avoid $t_j$.

For $\ell\ge26$, repeat the construction of Sections~\ref{sec:blocks}
and~\ref{sec:proof} with $G'$ in place of $G$ and
$L=a_1a_2\dots a_{k-11}$ for $14\le k\le\ell$.
In~\ref{it:F1} use $a_{k-11}$, and in~\ref{it:F3} replace $\{p,p'\}$ by $\{p^*\}$.
Each block has $13+(k-11)+(\ell-k)+1=\ell+3$ vertices.
The proof of Lemma~\ref{lem:block} needs only the following changes.
In part~\ref{it:block-all}, a path which enters $G'$ from $L$ has at
most $11$ vertices in $G'$, which gives the bound
$(k-11)+11+(\ell-k)=\ell$, and a path which does not meet $L$ has
order at most $13+(\ell-k)\le\ell-1$.
In part~\ref{it:block-in}, a path contained in $G'$ has order at most
$13\le k-1$, and a path which enters $G'$ from $L$ and ends in $G'$
outside $T$ has at most $10$ vertices in $G'$, which gives the bound
$(k-11)+10=k-1$.
In part~\ref{it:block-wit}, the paths of order $11$ obtained from the
paths $Q_j$ are used.
Since $k_i\ge(\ell+1)/2>13$, all blocks satisfy $k_i\ge14$.
The resulting digraph has $\ell(\ell+3)$ vertices and longest paths of order
$\ell$, and no stable set meets every longest path. For $\ell=26$ it has $754$ vertices.
\end{remark}

The oriented graphs $D_\ell$ disprove Conjecture~\ref{conj:DPPC} for $q=1$ and, by exchanging the two parts, for $q=\ell-1$.
This leaves the following question.

\begin{question}\label{q:middle}
Is it true that for every digraph $D$ and every integer $q$ with $2\le q\le\lambda(D)-2$, there exists a partition $(A,B)$ of $V(D)$ such that $\lambda(D[A])\le q$ and $\lambda(D[B])\le\lambda(D)-q$?
\end{question}

Conjecture~\ref{conj:LPX} asks for a set $A$ with $\lambda(D[A])\le1$ which meets every longest path of $D$. We may relax the bound on $\lambda(D[A])$ and ask the following question.

\begin{question}\label{q:constant}
Is there an absolute constant $c$ such that every digraph $D$ has a set $A\subseteq V(D)$ that meets every longest path of $D$ and satisfies $\lambda(D[A])\le c$?
\end{question}

By Theorem~\ref{thm:main}, any such constant $c$ must be at least $2$. A positive answer to Question~\ref{q:middle} would imply a positive
answer to Question~\ref{q:constant} with $c=2$. Indeed, if $\lambda(D)\ge4$,
take $q=2$. Then $\lambda(D[A])\le2$ and $\lambda(D[B])\le\lambda(D)-2$, so $A$ meets every longest path of $D$. If
$\lambda(D)\le3$, then we may take $A$ to be the stable set given by the following proposition. The problem page~\cite{REGS} records, without proof, that for such digraphs one can even require every vertex of the stable set to be the first vertex of a longest path. For completeness, we give a short proof of the proposition.

\begin{proposition}\label{prop:small}
Every digraph $D$ with $\lambda(D)\le3$ has a stable set which meets every longest path of $D$.
\end{proposition}

\begin{proof}
We may assume that $D$ is non-empty. A \emph{weak component} of $D$ is a subdigraph of $D$ induced by the vertex set of a connected component of the underlying graph of $D$. Suppose that $D$ has a directed triangle with arcs $v_1v_2$, $v_2v_3$ and $v_3v_1$, and let $v$ be a vertex of $D$ other than $v_1$, $v_2$ and $v_3$. If $vv_1$ were an arc of $D$, then $vv_1v_2v_3$ would be a path of order $4$, and if $v_1v$ were an arc of $D$, then $v_2v_3v_1v$ would be a path of order $4$. By symmetry, $v$ is adjacent to none of $v_1$, $v_2$ and $v_3$. Hence the three vertices of every directed triangle of $D$ form the vertex set of a weak component of $D$. We call such a weak component \emph{triangular}.

Let $D'$ be the digraph obtained from $D$ by deleting the vertices of all triangular weak components. Since $\lambda(D)\le3$, every directed cycle of $D$ has length $2$ or $3$, and so $D'$ has no directed cycle of odd length. If $D'$ is non-empty, then it has a kernel $N$ by Richardson's theorem~\cite{Richardson}; otherwise let $N=\emptyset$. For each triangular weak component $C$ of $D$, choose a vertex $s_C\in V(C)$, and let $S$ consist of the vertices of $N$ and the chosen vertices $s_C$. Since there is no arc between distinct weak components, $S$ is stable.

Let $P$ be a longest path of $D$. Then $P$ is contained in a weak component of $D$. If this weak component is triangular, then it has exactly three vertices and contains a path of order $3$, so $|P|=\lambda(D)=3$ and $P$ contains all three of its vertices; in particular, $P$ meets $S$. Otherwise $P$ is a path of $D'$ with $|P|=\lambda(D)\ge\lambda(D')$, so $P$ is a longest path of $D'$, and $N$ meets $P$ by the observation on kernels made in the Introduction.
\end{proof}

Finally, we ask how short the longest paths in a counterexample to
Conjecture~\ref{conj:LPX} can be. For oriented graphs, the minimum
possible value of $\lambda(D)$ lies between $6$ and $28$.
The lower bound follows from~\cite[Theorem~4.2]{GGM}, which covers
oriented graphs whose longest paths have at most four arcs (five vertices).
The upper bound follows from Theorem~\ref{thm:main}. For general digraphs, Proposition~\ref{prop:small} gives a lower bound of $4$, and Remark~\ref{rem:two-cycles} gives an upper bound of $26$.
\begin{question}\label{q:smallest} What is the smallest integer $\ell$ for which there exists a digraph (respectively, an oriented graph) $D$ with $\lambda(D)=\ell$ such that no stable set of $D$ meets every longest path of $D$?
\end{question}

\section*{Declaration of generative AI and AI-assisted technologies}

The authors proposed the construction strategy. Following this strategy,
OpenAI's GPT-6.2 Sol found the initial counterexample. It is the digraph
in Remark~\ref{rem:two-cycles} with $\ell=27$ and has $810$ vertices
and $81$ directed $2$-cycles. GPT-6.2 Sol also drew the diagrams in
Figures~\ref{fig:gadget}--\ref{fig:global}.

The authors verified the mathematical arguments, checked the computational
results and figures, revised the manuscript, and take full responsibility
for the content of the publication.

\end{document}